\documentclass[12pt]{article}
\usepackage[english]{babel}
\usepackage{amssymb}
\usepackage{amsmath}
\usepackage{graphicx}
\usepackage{bbm}
\usepackage[right]{eurosym}

\def\vs{\vspace}
\def\noi{\noindent}

\def\IN{\mathbb N}

\def\IR{\mathbb R}

\def\IQ{\mathbb Q}

\def\an{\mathrm{an}}
\def\exp{\mathrm{exp}}

\def\ma{\mathcal}

\begin{document}
	\begin{center}
		{\bf \Large The Derivative of a  
			Constructible Function 
			
			\vs{0.2cm}
			
			is Constructible}
	\end{center}
	
	\centerline{Tobias Kaiser}
	
	\vspace{0.7cm}\noi \footnotesize {{\bf Abstract.}
		The notion of constructible functions in the setting of tame real geometry has been introduced by Cluckers and Dan Miller in their work on parametric integration of globally subanalytic functions. A function on a globally subanalytic set is called constructible if it is a finite sum of finite products of globally subanalytic functions and the logarithm of positive globally subanalytic functions. We show that the class of constructible functions is stable under taking derivatives. 
		
		\normalsize
		\section*{Introduction}
		Log-analytic functions have been defined by Lion and Rolin in their seminal paper [10]. 
		They are iterated compositions from either side of globally subanalytic functions and the global logarithm and allow a preparation theorem.
		
		An important subclass is given by the so-called constructible functions. This notion has been coined by Cluckers and Dan Miller in their work on parametric integration of globally subanalytic functions [2, 3, 4], succeeding the work of Lion, Rolin [12] and Comte, Lion and Rolin [5]. A function on a globally subanalytic set is called constructible if it is a finite sum of finite products of globally subanalytic functions and the logarithm of positive globally subanalytic functions. The key property of the class of constructible functions is its closedness under parametric integration. For establishing this, Cluckers and Dan Miller have also proved a preparation theorem and a decay result in this setting.
		
		In [10] it was shown that the derivative of a differentiable log-analytic function  is again log-analytic.  In this note we show that the subclass of constructible functions is also stable under taking derivatives. 
		Here is the main result.

		\vs{0.5cm}
		{\bf Theorem}
		
		\vs{0.1cm}
		{\it Let $U\subset \IR^n$ be an open and globally subanalytic set and let $f:U\to \IR$ be constructible. Let $i\in \{1,\ldots,n\}$ be such that $f$ is differentiable with respect to the variable $x_i$ on $U$.  Then $\partial f/\partial x_i$ is constructible.}
		
			\vs{0.5cm}
		\hrule
		
		\vs{0.4cm}
		{\footnotesize{\itshape 2020 Mathematics Subject Classification:} 03C64, 14P15, 26A09, 32B20}
		\newline
		{\footnotesize{\itshape Keywords and phrases:} constructible functions, preparation, differentiability}
		
		\newpage
		This result is a necessary ingredient in the recent work of Huber, Oswal and the author [9], establishing a constructible de Rham theorem in the globally subanalytic setting.
		The proof of the theorem here adapts the strategy of [10, Theorem A], establishing a suitable preparation theorem for constructible functions on so-called simple cells. The preparation theorem is an appropriate variant of the one formulated in [3, Corollary 3.5]. Aizenbud et al. follow in their preprint [1, Propositions 2.7 \& 4.12] a similar approach to show that a much larger class called $\ma{C}^\exp$ is stable under taking limits and derivatives.

		\section*{Notations and Preliminaries}
		
		The empty sum is by definition $0$ and the empty product is by definition $1$.
		
		By $\IN=\{1,2,\ldots\}$ we denote the set of natural numbers and by $\IN_0=\{0,1,2,\ldots\}$ the set of nonnegative integers.

		We set $\IR_{>0}:=\{x\in \IR\mid x>0\}$. 
		For $a,b\in \IR$ with $a\leq b$ we denote by $[a,b]$ the closed interval and by $(a,b)$ the open interval with endpoints $a,b$, respectively.
		By $|\;\;|$ we denote the euclidean norm on $\IR^n$
		
		\vs{0.5cm}
		We assume familiarity with the notion of o-minimal structures and of globally subanalytic sets and functions (see for example van den Dries [6] or van den Dries and Miller [7]; see also the preliminaries of [10]). 
		
		\vs{0.5cm}
		We give the definition of a constructible function  established by Cluckers and D. Miller  in [2].
		
		\vs{0.5cm}
		{\bf Definition} 
		
		\vs{0.1cm}
		Let $X \subset \IR^n$ be a globally subanalytic set. A function $f:X \to \IR$ is called {\bf constructible} (on $X$) if there are $k,l\in \IN_0$ and globally subanalytic functions $(f_i:X\to \IR)_{1\leq i\leq k}$ and $(g_{ij}:X\to \IR_{>0})_{1\leq i\leq k, 1\leq j\leq l}$ such that 
		$$f=\sum_{i=1}^k f_i\prod_{j=1}^l \log(g_{ij}).$$
		
		\vs{0.5cm}
		Note that a constructible function is definable in the structure $\IR_{\an,\exp}$ and that the set of constructible functions on a given globally subanalytic set is a ring with respect to pointwise addition and multiplication which contains the globally subanalytic functions on that set. Note also that the class of constructible functions is closed under plugging in globally subanalytic maps.

		\section*{Results}

		We let $x=(x_1,\ldots,x_n)$ range over $\IR^n$ and $y$ over $\IR$. 
		We set $\pi:\IR^n\times \IR\to \IR^n, (x,y)\mapsto x$.
		
		\vs{0.5cm}
		For the following definition compare with [11, 2] and also with van den Dries and Speissegger [8]. Let $C\subset \IR^n\times\IR$ be globally subanalytic.

		\vs{0.5cm}
		{\bf 1. Definition}
		
		\vs{0.1cm}
		A function $\widetilde{y}:C\to \IR$ is called a {\bf globally subanalytic $y$-coordinate} on $C$ with {\bf center} $\Theta$ if the following holds:
		\begin{itemize}
			\item[(a)] We have $\widetilde{y}>0$ or $\widetilde{y}<0$ on $C$.
			\item[(b)] The function $\Theta$ is a globally subanalytic function on $\pi(C)$. 
			\item[(c)] We have $\widetilde{y}(x,y)=y-\Theta(x)$.
			\item[(d)] There is $\varepsilon \in (0,1)$ such that $0<|\widetilde{y}(x,y)|<\varepsilon|y|$ for all $(x,y) \in C$ or $\Theta = 0$.
		\end{itemize}
		
		\vs{0.2cm}
		Note that $\Theta$ is uniquely determined by $\widetilde{y}$.
		
		\vs{0.5cm}
		{\bf 2. Remark}
		
		\vs{0.1cm}
		Nota that Cluckers and D. Miller omit condition (d) in their results [3, Theorem 3.4 \& Corollary 3.5]. But we can add it without problems by the original formulation of the globally subanalytic preparation theorem in [11, p. 862] or by [8, Theorem 2.4].
		
		\vs{0.5cm}
		We give the notion of so-called strong functions from (cf. [3, Definition 3.3]).

		\vs{0.5cm}
		{\bf 3. Definition}
		
		\vs{0.1cm}
		A function $u:C \to \IR$ is called  {\bf strong} on $C$ with globally subanalytic $y$-coordinate $\widetilde{y}$ on $C$ if $u=v \circ \varphi$ where the following holds:
		\begin{itemize}
			\item [(a)] The function $\varphi$ is given by
			\begin{eqnarray*}
				\varphi:C &\to& [-1,1]^s,  \\
				\varphi(x,y)	&=&\Big(b_1(x),\ldots,b_{s-2}(x),b_{s-1}(x)|\widetilde{y}|^p,b_s(x)|\widetilde{y}|^{-p}\Big),
			\end{eqnarray*}
			where $s\in \IN_0$, $p \in \IQ$ and $b_1,...,b_s$ are globally subanalytic functions on $\pi(C)$.
			\item[(b)] The function $v$ is a power series which converges on a neighbourhood of $[-1,1]^s$.
		\end{itemize}
		
		\vs{0.2cm}
		Note that a strong function is globally subanalytic and bounded.
		
		\vs{0.5cm}
		Given $x\in \IR^n$, we set $C_x:=\{y\in \IR\mid (x,y)\in C\}$. We introduce the notion of a simple cell (see [10, Definition 2.15]). 
		
		\vs{0.5cm}
		{\bf 4. Definition}
		
		\vs{0.1cm}
		We call $C$ {\bf simple} if for every $x\in \pi(C)$ we have $C_x=(0,d_x)$ for some $0<d_x<+\infty$.
		
		\vs{0.5cm}
		{\bf 5. Remark}
		
		\vs{0.1cm}
		Let $X\subset \IR^n$ be a globally subanalytic set and 
		let $\ma{C}$ be a globally subanalytic cell decomposition of $X\times (0,1)$. Then
		$$X=\bigcup\{\pi(C)\mid C\in\ma{C}\mbox{ simple}\}.$$ 
		
		\vs{0.5cm}
		{\bf 6. Remark}
		
		\vs{0.1cm} 
		Let $C$ be simple. Then $y$ is the only globally subanalytic $y$-coordinate on $C$. 
		
		\vs{0.1cm}
		{\bf Proof:}
		
		\vs{0.1cm}
		That $y$ is a globally subanalytic $y$-coordinate on $C$ is clear.
		Let $\widetilde{y}$ be a globally subanalytic $y$-coordinate on $C$ with center $\Theta$.    
		Assume that $\Theta \neq 0$. Then by Definition 3 there is $\varepsilon \in (0,1)$ such that
		$$|y-\Theta(x)| < \varepsilon|y|$$
		for all $(x,y) \in C$. Let $x \in \pi(C)$ be such that $\Theta(x) \neq 0$. Then we obtain
		$$+\infty=\lim_{y\searrow 0}\Big\vert 1-\frac{\Theta(x)}{y}\Big\vert\leq \varepsilon$$ 
		which is a contradiction.  
		\hfill$\blacksquare$
		
		\vs{0.5cm}
		{\bf 7. Theorem}
		
		\vs{0.1cm}
		{\it Let $X\subset \IR^n$ be a globally subanalytic set and let $F:X\times (0,1)$ be a constructible function. 
			Then there exists a cell decomposition $\ma{C}$ of $X\times (0,1)$ such that for every $C\in \ma{C}$ which is simple the following holds:
			\begin{itemize}
				\item[(1)] We have 
				$$f|_C(x,y)=\sum_{i=1}^M a_i(x)u_i(x,y) y^{p_i} |\log y|^{l_i}$$
				where $M\in \IN_0$ and for each $i\in \{1,\ldots,M\}$, $a_i:\pi(C)\to \IR$ is constructible, $u_i:C\to \IR$ is strong, $p_i\in \IQ$ and $l_i\in \IN_0$.
				\item[(2)] For each $i\in \{1,\ldots,M\}$ we have $u_i=1$ or $p_i>0$.
				\item[(3)] The pairs $(p_i,l_i)$ where $p_i\leq 0$ are pairwise distinct.
		\end{itemize}}
		
		\vs{0.1cm}
		{\bf Proof:}
		
		\vs{0.1cm}
		By [3, Corollary 3.5] in connection with Remark 3 and Remark 6 we get such a cell decomposition such that Property (1) holds.
		Applying the additional statment of [3, Corollary 3.5] we immediately obtain a presentation such that for each $i$ we have $u_i=1$ or $p_i>-1$. But the proof thereof can be exactly used to obtain also Property (2). Property (3) is then in view of Property (2) easily obtained by collecting summands.
		\hfill$\blacksquare$

		\vs{0.5cm}
		{\bf 8. Theorem}
		
		\vs{0.1cm}
		{\it Let $X \subset \IR^n$ be globally subanalytic and let $F:X\times (0,1)\to \IR$ be constructible. Assume that for every $x\in X$ we have that $\lim_{y\searrow 0}F(x,y)$ exists and is finite.
			Then the function $h:X\to \IR, x\mapsto \lim_{y\searrow 0}F(x,y),$ is constructible.}
		
		\vs{0.1cm}
		{\bf Proof:}
		
		\vs{0.1cm}
		By Theorem 7 we find a cell decomposition $\ma{C}$ of $X\times (0,1)$ such that for every simple $C\in\ma{C}$ 
		we have 
		$$F|_C(x,y)=\sum_{i=1}^M a_i(x)u_i(x,y) y^{p_i} |\log y|^{l_i}.$$
		with the properties mentioned therein.
		We show that $h|_{\pi(C)}$ is constructible for such $C$ and are done by Remark 5.
		For $x\in \pi(C)$ and $i\in \{1,\ldots,M\}$ with $p_i>0$ we have
		$$\lim_{y\searrow 0} a_i(x)u_i(x,y) |y|^{p_i}|\log y|^{l_i}=0.$$
		Hence we can replace $F|_C$ by 
		$$\sum_{p_i\leq 0} a_i(x)u_i(x,y) y^{p_i} |\log y|^{l_i}=\sum_{p_i\leq 0} a_i(x) y^{p_i} |\log y|^{l_i}.$$
		We can therefore assume that a priori $p_i\leq 0$ for all $i\in \{1,\ldots,M\}$. We may also assume that $M\geq 1$ and that $a_i$ does not vanish identically for each $i\in \{1,\ldots,M\}$.
		We may further assume that $p_1$ is minimal among the $p_i$ and that $l_1$ is maximal among those $l_i$ for which $p_i=p_1$.
		Then 
		$$\lim_{y\searrow 0} \frac{y^{p_i}|\log y|^{l_i}}{y^{p_1}|\log y|^{l_1}}=0$$
		for every $i>1$. Since $a_1$ does not vanish identically we obtain by the assumption $\lim_{y\searrow 0} F(x,y)\in \IR$ for every $x\in \pi(C)$ that $(p_1,l_1)=(0,0)$ (and that $M=1$).
		We conclude that $h|_{\pi(C)}=a_1$ and are done.
		\hfill$\blacksquare$

		\vs{0.5cm}
		With the above theorem we are able to establish our main result (compare the proof of [10, Theorem A]).
		
		\vs{0.5cm}
		{\bf 9. Theorem}
		
		\vs{0.1cm}
		{\it Let $U\subset \IR^n$ be an open and globally subanalytic set and let $f:U\to \IR$ be constructible. Let $i\in \{1,\ldots,n\}$ be such that $f$ is differentiable with respect to the variable $x_i$ on $U$.  Then $\partial f/\partial x_i$ is constructible.}
		
		\vs{0.1cm}
		{\bf Proof:}
		
		\vs{0.1cm}
		We may assume that $f$ is differentiable with respect to the last variable $x_n$. We have to show that $\partial f/\partial x_n$ is constructible.
		Let $e_n:=(0,\ldots,0,1)\in \IR^n$ be the $n^\mathrm{th}$ unit vector.
		We define
		$$F:U\times (0,1)\to \IR, (x,y)\mapsto \left\{\begin{array}{ccc}
			\frac{f(x+ye_n)-f(x)}{y},&& x+ye_n\in U,\\
			&\mbox{if}&\\
			0,&& x+ye_n\notin U.
		\end{array}	\right.$$
		Then $F$ is constructible. 
		Since
		$$\frac{\partial f}{\partial x_n}(x)=\lim_{y\searrow 0}F(x,y)$$
		for $x\in U$ we are done by Theorem 8.
		\hfill$\blacksquare$

		\vs{1cm}
		\noi \footnotesize{\centerline{\bf References}
			\begin{itemize}
				\item[(1)] A. Aizenbud, R. Cluckers, M. Raibaut, T. Servi:
				Analytic holonomicity of real $\ma{C}\exp$-class distributions.
				arXiv 2403.20176.
				\item[(2)] 
				R. Cluckers and D. Miller:
				Stability under integration of sums of products of real globally subanalytic functions and their logarithms.
				{\it Duke Math. J.} {\bf 156} (2011), no. 2, 311-348.
				\item[(3)] 
				R. Cluckers and D. Miller:
				Loci of integrability, zero loci, and stability under integration for constructible functions on Euclidean space with Lebesgue measure.
				{\it Int. Math. Res. Not.} (2012), no. 14, 3182-3191.
				\item[(4)]
				R. Cluckers and D. Miller:
				Lebesgue classes and preparation of real constructible functions.
				{\it J. Funct. Anal.} {\bf 264} (2013), no. 7, 1599-1642.
				\item[(5)] 
				G. Comte, J.-M. Lion and J.-P. Rolin:
				Nature log-analytique du volume des sous-analytiques.
				{\it Illinois J. Math.} {\bf 44} (2000), no. 4, 884-888.
				\item[(6)]
				L. van den Dries: Tame Topology and O-minimal Structures. {\it London Math. Soc. Lecture Notes Series} {\bf 248}, Cambridge University Press, 1998.
				\item[(7)]
				L. van den Dries and C. Miller:
				Geometric categories and o-minimal structures.
				{\it Duke Math. J.} {\bf 84} (1996), no. 2, 497-540.
				\item[(8)]
				L. van den Dries and P. Speissegger:
				O-minimal preparation theorems.
				Model theory and applications, 87-116, Quad. Mat., 11, Aracne, Rome, 2002.
				\item[(9)] 
				A. Huber, T. Kaiser, A. Oswal: On the de Rham Theorem in the Globally Subanalytic Setting. 
				arXiv 2508.03499, 27 p.
				\item[(10)]
				T. Kaiser, A. Opris: Differentiability Properties of Log-Analytic Functions.
				{\it Rocky Mountain Journal of Mathematics} {\bf 52} (2022) no. 4, 1423-1443.
				\item[(11)]
				J.-M. Lion, J.-P. Rolin:
				Th\'{e}or\`{e}me de pr\'{e}paration pour les fonctions logarithmico-exponentielles.
				{\it Ann. Inst. Fourier} {\bf 47}, no. 3 (1997), 859-884.
				\item[(12)]
				J.-M. Lion, J.-P. Rolin:
				Int\'egration des fonctions sous-analytices et volumes des sous-ensembles sous-analytiques. 
				{\it Ann. Inst. Fourier} {\bf 48}, no. 3 (1998), 755-767.
		\end{itemize}}
		
		\vs{0.5cm}
		Tobias Kaiser\\
		University of Passau\\
		Faculty of Computer Science and Mathematics\\
		tobias.kaiser@uni-passau.de\\
		D-94030 Germany

	\end{document}